\input amstex
\input amsppt.sty
\magnification=\magstep1
\hsize=30truecc
\vsize=22.2truecm
\baselineskip=16truept
\TagsOnRight
\pageno=1
\nologo
\def\Z{\Bbb Z}
\def\N{\Bbb N}

\def\Q{\Bbb Q}

\def\l{\left}
\def\r{\right}
\def\bg{\bigg}
\def\({\bg(}
\def\[{\bg\lfloor}
\def\){\bg)}
\def\]{\bg\rfloor}
\def\t{\text}
\def\f{\frac}

\def\bi{\binom}
\def\eq{\equiv}

\def\ls{\leqslant}
\def\gs{\geqslant}
\def\mo{\roman{mod}}
\def\ord{\roman{ord}}

\def\al{\alpha}
\def\da{\delta}

\def\Proof{\noindent{\it Proof}}

\def\Ack{\medskip\noindent {\bf Acknowledgments}}
\hbox {Acta Arith. 122(2006), no.\,1, 91--100.}
\bigskip
\topmatter
\title Polynomial extension of Fleck's congruence\endtitle
\author Zhi-Wei Sun (Nanjing)\endauthor
\leftheadtext{Zhi-Wei Sun}
\abstract
Let $p$ be a prime, and let $f(x)$ be an integer-valued polynomial.
By a combinatorial approach, we obtain a nontrivial lower bound of
the $p$-adic order of the sum
$$\sum_{k\eq r\, (\mo\ p^{\beta})}\bi nk(-1)^kf\l(\l\lfloor\f{k-r}{p^{\al}}\r\rfloor\r),$$
where $\al\gs\beta\gs0$, $n\gs p^{\al-1}$ and $r\in\Z$.
This polynomial extension
of Fleck's congruence has various backgrounds and several consequences
such as
$$\sum_{k\eq r\,(\mo\ p^\al)}\bi nka^k\eq0\
\(\mo\ p^{\big\lfloor\f{n-p^{\al-1}}{\varphi(p^\al)}\big\rfloor}\)$$
provided that $\al>1$ and $a\eq-1\ (\mo\ p)$.
\endabstract
\thanks 2000 {\it Mathematics Subject Classification}:\,Primary 11B65;
Secondary 05A10, 11A07, 11B68, 11S05.\newline\indent
The author is partially supported by the National Science Fund
for Distinguished Young Scholars (No. 10425103) and a Key
Program of NSF in P. R. China.
\endthanks
\endtopmatter
\document

\heading{1. Introduction}\endheading

As usual, we let $\bi x0=1$ and
$$\bi xk=\f{x(x-1)\cdots(x-k+1)}{k!}\quad \t{for every}\ k=1,2,3,\ldots.$$
For convenience, we also set $\bi xk=0$ for any negative integer $k$.

 Let $p$ be a prime and $r$ be an integer.
 In 1913 A. Fleck (cf. Dickson [D, p.\,274]) discovered that
 $$\sum_{k\eq r\, (\mo\ p)}\bi nk(-1)^k\eq0
 \ \l(\mo\ p^{\l\lfloor\f{n-1}{p-1}\r\rfloor}\r)\tag1.1$$
 for all $n\in\Z^+=\{1,2,3,\ldots\}$,
 where $\lfloor\cdot\rfloor$ is the well-known floor function.
 Sums of the form $\sum_{k\eq r\, (\mo\ m)}\bi nk$ or $\sum_{k\eq r\, (\mo\ m)}\bi nk(-1)^k$
 (with $m\in\Z^+$)
 have various applications in number theory and combinatorics
 (see, e.g., [SS], [H] and [S02]).

 In 1977, by a very complicated method, C. S. Weisman [W]
extended Fleck's congruence to prime power moduli in the following way:
$$\sum_{k\eq r\, (\mo\ p^{\al})}\bi nk(-1)^k
\eq0\ \l(\mo\ p^{\big\lfloor\f{n-p^{\al-1}}{\varphi(p^{\al})}\big\rfloor}\r),\tag1.2$$
where $\al,n\in\N=\{0,1,2,\ldots\}$ and $n\gs p^{\al-1}$, and $\varphi$
denotes Euler's totient function.
Unaware of Fleck's previous work, Weisman was motivated by studying the relation
between two different ways (Mahler's and van der Put's) to express
a $p$-adically continuous function.

 Quite recently,
in his lecture notes on Fontaine's rings and $p$-adic
$L$-functions given at Irvine (Spring, 2005), D. Wan got the following
new extension of Fleck's congruence:
$$\sum_{k\eq r\, (\mo\ p)}\bi nk(-1)^k\bi{(k-r)/p}l
\eq0\ \l(\mo\ p^{\l\lfloor\f{n-lp-1}{p-1}\r\rfloor}\r),\tag1.3$$
where $l,n\in\N$ and $n>lp$.
Wan was led to this when trying to understand a sharp estimate for the $\psi$-operator
in Fontaine's theory of ($\phi$,\,$\Gamma$)-modules.

 For a prime $p$, we let $\Q_p$ and $\Z_p$ denote the field of $p$-adic numbers
 and the ring of $p$-adic integers respectively;
the $p$-adic order of $\omega\in\Q_p$ is defined by
$\ord_p(\omega)=\sup\{a\in\Z:\, \omega/p^{a}\in\Z_p\}$ (whence $\ord_p(0)=+\infty$).
Throughout this paper,
the Kronecker symbol $\da_{m,n}$ with $m,n\in\N$ takes $1$ or $0$
according as $m=n$ or not.

 Clearly both Weisman's and Wan's extensions of Fleck's congruence
 follow from the special case $\al=\beta$ of the following theorem,
 which we will establish by a combinatorial approach.

\proclaim{Theorem 1.1} Let $p$ be a prime, and let $f(x)\in\Q_p[x]$,
$\deg f\ls l\in\N$ and $f(a)\in\Z_p$ for all $a\in\Z$.
Provided that $\al,\beta\in\N$ and $\al\gs\beta$,
we have
$$\sum_{k\eq r\, (\mo\ p^{\beta})}\bi nk
(-1)^kf\l(\l\lfloor\f{k-r}{p^{\al}}\r\rfloor\r)\in
p^{\big\lfloor\f{n-p^{\al-1}-l}{\varphi(p^{\al})}\big\rfloor-(l-1)\al-\beta}\Z_p\tag1.4$$
for all integers $n\gs p^{\al-1}$ and $r$;
moreover, we can substitute $\da_{\beta,0}$ for the first $l$
in $(1.4)$ if $\al$ is greater than one.
\endproclaim

By Theorem 1.1 in the case $\al=\beta=r=0$, if $f(x)\in\Z[x]$ and $f(x)\not=0$, then
for any integer $n>\deg f+1$ we have $\sum_{k=0}^n\bi nk(-1)^kf(k)=0$
since the sum is divisible by all primes. In fact, a known identity due to L. Euler
(cf. [LW, pp.\,90-91]) states that
$$\sum_{k=0}^n\bi nk(-1)^{k}k^l=\cases(-1)^n\, n!&\t{if}\ l=n\in\N,
\\0&\t{if}\ 0\ls l<n.\endcases$$

Now we derive more consequences of Theorem 1.1.

\proclaim{Corollary 1.1} Let $p$ be a prime, $m\in\Z^+$ and $\al=\ord_p(m)$.
Let $l,n\in\N$ and $r\in\Z$. Then
$$\aligned&\ord_p\(\sum_{k\eq r\,(\mo\ p^{\al})}\bi nk(-1)^kB_l\l(\f{k-r}m\r)\)
\\&\quad\gs\l\lfloor\f{n-p^{\al-1}-l(\da_{\al,0}+\da_{\al,1})}
{\varphi(p^{\al})}\r\rfloor-l\al,
\endaligned\tag1.5$$
where $B_l(x)$ is the Bernoulli polynomial of degree $l$.
\endproclaim
\Proof. (1.5) holds trivially if $n<p^{\al-1}$. Below we suppose $n\gs p^{\al-1}$.

When $l=0$, (1.5) reduces to Weisman's congruence (1.2).
In the case $\al=0$, if the lower bound in (1.5) is nonnegative (i.e., $l<n$)
then the summation in (1.5)
vanishes by Euler's identity.

Now we assume $l\al\not=0$, and let
$B_l=B_l(0)$ be the $l$th Bernoulli number.
Note that $m_0=m/p^{\al}$ is relatively prime to $p$.
For any $a\in\Z$ we have $B_l(a/m_0)-B_l\in\Z_p$, because
$$m_0^l\l(B_l\l(\f a{m_0}\r)-B_l\r)
=\l(m_0^lB_l\l(\f a{m_0}\r)-B_l\r)-\l(m_0^lB_l(0)-B_l\r)\in\Z_p$$
by [S03, Corollary 1.3].
Applying Theorem 1.1 with $f(x)=B_l(x/m_0)-B_l$ and $\beta=\al$, we get that
$$\align&\ord_p\(\sum_{k\eq r\,(\mo\ p^{\al})}\bi nk(-1)^kB_l\l(\f{k-r}m\r)-B_l\Sigma\)
\\&\qquad\quad
\gs\l\lfloor\f{n-p^{\al-1}-l\da_{\al,1}}{\varphi(p^{\al})}\r\rfloor-l\al,
\endalign$$
where $\Sigma=\sum_{k\eq r\,(\mo\ p^{\al})}\bi nk(-1)^k$.
Recall that $pB_l\in\Z_p$ by the von Staudt--Clausen theorem (cf. [IR, pp.\,233-236]).
This, together with (1.2), shows that
$$\ord_p(B_l\Sigma)\gs\ord_p(\Sigma)-1
\gs\l\lfloor\f{n-p^{\al-1}}{\varphi(p^{\al})}\r\rfloor-1
\gs\l\lfloor\f{n-p^{\al-1}-l\da_{\al,1}}{\varphi(p^{\al})}\r\rfloor-l\al.$$
So the desired (1.5) follows.
\qed

\proclaim{Corollary 1.2} Let $p$ be a prime, and let $f(x)\in\Q_p[x],\ \deg f=l\gs0$
and $f(a)\in\Z_p$ for all $a\in\Z$. Let $\al\in\N$ and $r\in\Z$.
Then, for any integer $n\gs p^{\al-1}$, we have
$$\aligned&\ord_p\(\sum_{k=0}^n\bi nk(-1)^k
(k-r,p^{\al})f\l(\l\lfloor\f{k-r}{p^{\al}}\r\rfloor\r)\)
\\&\quad\gs\l\lfloor\f{n-p^{\al-1}-l(\da_{\al,0}+\da_{\al,1})}
{\varphi(p^{\al})}\r\rfloor-(l-1)\al-1,
\endaligned\tag1.6$$
where $(k-r,p^{\al})$ is the greatest common divisor of $k-r$ and $p^{\al}$.
\endproclaim
\Proof. Let $g(1)=p$ and $g(p^{\beta})=p-1$ if $0<\beta\ls\al$.
By Theorem 1.1, the $p$-adic order of
$$\align&\sum_{\beta=0}^{\al}g(p^{\beta})p^{\beta}\sum_{k\eq r\,(\mo\ p^{\beta})}\bi nk(-1)^k
f\l(\l\lfloor\f{k-r}{p^{\al}}\r\rfloor\r)
\\&=\sum_{k=0}^n\bi nk(-1)^kf\l(\l\lfloor\f{k-r}{p^{\al}}\r\rfloor\r)
\sum_{d\mid(k-r,p^{\al})}g(d)d
\endalign$$
is at least
$$\nu=\l\lfloor\f{n-p^{\al-1}-l(\da_{\al,0}+\da_{\al,1})}{\varphi(p^{\al})}\r\rfloor-(l-1)\al.$$
We note in passing that in the case $\al>1$,
$$\ord_p(g(p^0))+\l\lfloor\f{n-p^{\al-1}-\da_{0,0}}
{\varphi(p^\al)}\r\rfloor\gs\l\lfloor\f{n-p^{\al-1}}{\varphi(p^\al)}\r\rfloor.$$
Now, since
$$\sum_{d\mid(k-r,p^{\al})}g(d)d
=p+\sum_{1<d\mid(k-r,p^{\al})}(p-1)d
=\sum_{d\mid(k-r,p^{\al})}\varphi(d)p=(k-r,p^{\al})p,$$
by the above the sum in (1.6) has $p$-adic order at least $\nu-1$. \qed

\proclaim{Corollary 1.3} Let $p$ be a prime, and let $\al,\beta,a,n,r$
be integers for which
$$\al>1,\ \al\gs\beta\gs0,\ a\eq\ 1\ (\mo\ p^{\al}),\ n\gs p^{\al-1}\ \t{and}\ r<p^{\beta}.$$
Then we have the congruence
$$\sum_{k\eq r\,(\mo\ p^{\beta})}\bi nk(-1)^ka^{\lfloor\f {k-r}{p^{\al}}\rfloor}
\eq0\ \(\mo\ p^{\big\lfloor\f{n-p^{\al-1}-\da_{\beta,0}}
{\varphi(p^{\al})}\big\rfloor+\al-\beta}\).
\tag1.7$$
\endproclaim
\Proof. When $a=1$, (1.7) holds
by Theorem 1.1 in the case $l=0$. So it suffices to show that
$$D:=\sum_{k\eq r\,(\mo\ p^{\beta})}\bi nk(-1)^k\l(a^{\lfloor\f{k-r}{p^{\al}}\rfloor}-1\r)$$
is divisible by $p^{\lambda}$
where
$$\lambda=\l\lfloor\f{n-p^{\al-1}-\da_{\beta,0}}{\varphi(p^{\al})}\r\rfloor+\al-\beta.$$

Write $a=1+p^{\al}b$ with $b\in\Z$. Then
$$\align D=&\sum_{k\eq r\,(\mo\ p^{\beta})}\bi nk(-1)^k
\sum_{0<l\ls\lfloor\f {k-r}{p^{\al}}\rfloor}\bi{\lfloor(k-r)/p^{\al}\rfloor}l(p^{\al}b)^l
\\=&\sum_{0<l\ls\lfloor\f {n-r}{p^{\al}}\rfloor}p^{l\al}b^l
\sum_{k\eq r\,(\mo\ p^{\beta})}\bi nk(-1)^k\bi{\lfloor(k-r)/p^{\al}\rfloor}l.
\endalign$$
For each $0<l\ls\lfloor(n-r)/p^{\al}\rfloor$,
applying Theorem 1.1 with $f(x)=\bi xl$ we find that
$$p^{l\al}\sum_{k\eq r\,(\mo\ p^{\beta})}\bi nk(-1)^k
\bi{\lfloor(k-r)/p^{\al}\rfloor}l\eq0\ \l(\mo\ p^{\lambda}\r).$$
Therefore $D\eq0\ (\mo\ p^{\lambda})$.
This concludes the proof. \qed
\medskip

Let $a\in\Z$ be congruent to $1$ modulo a prime $p$.
By induction, $a^{p^{\al}}\eq1\ (\mo\ p^{\al+1})$ for any $\al\in\N$.
Let $n,r\in\Z$ and $n\gs p^{\al-1}$. If $\al\gs2$, then
by Corollary 1.3 in the case $\beta=\al$ we have
$$\sum_{k\eq r\,(\mo\ p^\al)}\bi nk(-a)^k\eq0\
\l(\mo\ p^{\big\lfloor\f{n-p^{\al-1}}{\varphi(p^\al)}\big\rfloor}\r).\tag1.8$$
By the binomial theorem, (1.8) is also valid with $\al=0$.
We remark that (1.8) also holds when $\al=1$, as pointed out by Fleck (cf. [D, p.\,274]).
\medskip

In the next section we will provide some lemmas.
Section 3 is devoted to our proof of Theorem 1.1.

\heading{2. Some Lemmas}\endheading

Let us recall the following well-known convolution identity of Chu and Vandermonde
(see, e.g., [GKP, (5.27)]):
$$\sum_{k=0}^n\bi xk\bi y{n-k}=\bi{x+y}n\quad \t{for all}\ n=0,1,2,\ldots.$$
This can be seen by
comparing the power series expansions of $(1+t)^x(1+t)^y$ and $(1+t)^{x+y}$.

\proclaim{Lemma 2.1} Let
$f(x)$ be a function from $\Z$ to a field, and let $m,n\in\Z^+$.
Then, for any $r\in\Z$ we have
$$\sum_{k=0}^n\bi nk(-1)^kf\l(\l\lfloor\f{k-r}m\r\rfloor\r)
=\sum_{k\eq \bar r\,(\mo\ m)}\bi{n-1}k(-1)^{k-1}\Delta f\l(\f{k-\bar r}m\r),$$
where $\bar r=r+m-1$ and $\Delta f(x)=f(x+1)-f(x)$.
\endproclaim
\Proof. By the Chu-Vandermonde identity, for any $h\in\N$ we have
$$\sum_{k=0}^h\bi nk(-1)^k=(-1)^h\sum_{k=0}^h\bi nk\bi{-1}{h-k}=(-1)^h\bi{n-1}h.$$
Therefore
$$\sum_{k=0}^n\bi nk(-1)^kf\l(\l\lfloor\f{k-r}m\r\rfloor\r)
=\sum_{j\in\Z}c_jf(j),$$
where
$$\align c_j=&\sum\Sb k\in\Z\\\lfloor\f{k-r}m\rfloor=j\endSb\bi nk(-1)^k
\\=&\sum_{0\ls k<(j+1)m+r}\bi nk(-1)^k-\sum_{0\ls k<jm+r}\bi nk(-1)^k
\\=&(-1)^{(j+1)m+r-1}\bi{n-1}{(j+1)m+r-1}
-(-1)^{jm+r-1}\bi{n-1}{jm+r-1}.
\endalign$$
(Note that $\bi{n-1}i\not=0$ only for those $i\in\{0,\ldots,n-1\}$.)
So we have
$$\align&\sum_{k=0}^n\bi nk(-1)^kf\l(\l\lfloor\f{k-r}m\r\rfloor\r)
\\=&\sum_{j\in\Z}(-1)^{(j+1)m+r-1}\bi{n-1}{(j+1)m+r-1}f(j)
\\&-\sum_{j\in\Z}
(-1)^{jm+r-1}\bi{n-1}{jm+r-1}f(j)
\\=&\sum_{k\eq \bar r\, (\mo\ m)}\bi{n-1}k(-1)^k\l(f\l(\f{k-\bar r}m\r)
-f\l(\f{k-\bar r}m+1\r)\r)
\\=&\sum_{k\eq \bar r\, (\mo\ m)}\bi{n-1}k(-1)^{k-1}\Delta f\l(\f{k-\bar r}m\r).
\endalign$$
This proves the desired identity. \qed
\medskip

 It is interesting to compare the identity in Lemma 2.1
with the following observation
$$\sum\Sb 0\ls k\ls n\\k\eq r\,(\mo\ m)\endSb\Delta f\l(\f{k-r}m\r)
=f\(\l\lfloor\f{n-r}m\r\rfloor+1\)-f\(\l\lfloor\f{-r-1}m\r\rfloor+1\),$$
which appeared in the author's proof of [S03, Lemma 3.1].

\proclaim{Lemma 2.2} Let $p$ be a prime and $\al$ be a positive integer. Then, for any
$k=0,1,\ldots,\varphi(p^{\al})$, we have
$$\bi{\varphi(p^{\al})}k\eq\cases(-1)^k\ (\mo\ p)&\t{if}\ p^{\al-1}\mid k,
\\0\ (\mo\ p)&\t{otherwise}.\endcases$$
\endproclaim
\Proof. Let $k=k_0+k_1p+\cdots+k_{\al-1}p^{\al-1}$ be the $p$-adic expansion
of $k$, where $k_0,k_1,\ldots,k_{\al-1}\in\{0,\ldots,p-1\}$. By a well-known theorem
of E. Lucas (see, e.g., [HS]),
$$\align\bi{\varphi(p^{\al})}k
=&\bi{\sum_{0\ls j<\al-1}0p^j+(p-1)p^{\al-1}}{\sum_{0\ls j<\al-1}k_jp^j+k_{\al-1}p^{\al-1}}
\\\eq&\bi{p-1}{k_{\al-1}}\prod_{0\ls j<\al-1}\bi 0{k_j}\ \ (\mo\ p).
\endalign$$

 If $p^{\al-1}\nmid k$, then $k_j>0$ for some $j<\al-1$, and hence
$\bi{\varphi(p^{\al})}k\eq0\ (\mo\ p)$. When $p^{\al-1}\mid k$, we have
$k_j=0$ for all $j<\al-1$, and thus
$$\align \bi{\varphi(p^{\al}}k&\eq\bi{p-1}{k_{\al-1}}
=\prod_{0<s\ls k_{\al-1}}\f{p-s}s\ \ (\mo\ p)
\\&\eq(-1)^{k_{\al-1}}\eq(-1)^{p^{\al-1}k_{\al-1}}=(-1)^k\ \ (\mo\ p).
\endalign$$
This completes the proof. \qed

\heading{3. Proof of Theorem 1.1}\endheading

We use induction on $w_l(\al,\beta):=l(\al+1)+\beta$.

In the case $w_l(\al,\beta)=0$ (i.e., $l=\beta=0$), the desired result is trivial because
$\sum_{k=0}^n\bi nk(-1)^k=(1-1)^n=0$ for all $n\in\Z^+$.

Let $w$ be a positive integer, and assume that the desired result holds
whenever $w_l(\al,\beta)<w$.
Now we deal with the case $w_l(\al,\beta)=w$.

{\bf Case 1}: $\beta=0$.

In this case, $l$ is positive.
Let $n\in\N$, $n\gs p^{\al-1}$,
$r\in\Z$ and $\bar r=r+p^{\al}-1$.
By Lemma 2.1,
$$\aligned&\sum_{k=0}^n\bi nk(-1)^kf\l(\l\lfloor\f{k-r}{p^{\al}}\r\rfloor\r)
\\=&\sum_{k\eq \bar r\,(\mo\ p^{\al})}\bi{n-1}k(-1)^{k-1}\Delta f\l(\f{k-\bar r}{p^{\al}}\r).
\endaligned\tag3.1$$
Clearly $\Delta f(x)$ is a
polynomial of degree at most $l-1$, and $\Delta f(a)\in\Z_p$ for all $a\in\Z$. Also,
$w_{l-1}(\al,\al)<w_l(\al,0)=w$.
In view of (3.1) and the induction hypothesis,
$$\align&\ord_p\(\sum_{k=0}^n\bi nk(-1)^kf\l(\l\lfloor\f{k-r}{p^{\al}}\r\rfloor\r)\)
\\\gs&\l\lfloor\f{(n-1)-p^{\al-1}-(l-1)}{\varphi(p^{\al})}\r\rfloor-(l-2)\al-\al
\\=&\l\lfloor\f{n-p^{\al-1}-l}{\varphi(p^{\al})}\r\rfloor-(l-1)\al-0.
\endalign$$
(Note that this is trivial if $n-1<p^{\al-1}$.)
Similarly, when $\al>1$, by (3.1) and the induction hypothesis we have
$$\align&\ord_p\(\sum_{k=0}^n\bi nk(-1)^kf\l(\l\lfloor\f{k-r}{p^{\al}}\r\rfloor\r)\)
\\\gs&\l\lfloor\f{(n-1)-p^{\al-1}-\da_{\al,0}}{\varphi(p^{\al})}\r\rfloor-(l-2)\al-\al
\\=&\l\lfloor\f{n-p^{\al-1}-\da_{0,0}}{\varphi(p^{\al})}\r\rfloor-(l-1)\al-0.
\endalign$$

{\bf Case 2}: $0<\beta\ls\al$.

If $l=0$ (i.e., $f(x)$ is constant), then
$w_l(\beta,\beta)=w_l(\al,\beta)=w$ and it suffices to handle the case $\al=\beta$.
In fact, when $l=0$, $n\gs p^{\al-1}$ and $r\in\Z$, provided that
$$\sum_{k\eq r\,(\mo\ p^{\beta})}\bi nk(-1)^kf\l(\f{k-r}{p^\beta}\r)
\in p^{\big\lfloor\f{n-p^{\beta-1}}{\varphi(p^\beta)}\big\rfloor}\Z_p$$
we have
$$\sum_{k\eq r\,(\mo\ p^{\beta})}\bi nk(-1)^kf\l(\l\lfloor\f{k-r}{p^\al}\r\rfloor\r)
\in p^{\big\lfloor\f{n-p^{\al-1}}{\varphi(p^\al)}\big\rfloor-(0-1)\al-\beta}\Z_p,$$
because
$$\f{n-p^{\beta-1}}{\varphi(p^\beta)}-\f{n-p^{\al-1}}{\varphi(p^\al)}
=\f n{p^{\al-1}}\sum_{0\ls s<\al-\beta}p^s\gs\al-\beta.$$
Below we simply let $(l-1)\al+\beta\gs0$ (i.e., $\al=\beta$ if $l=0$).

Let us use induction on $n\gs p^{\al-1}$.
The desired result is trivial when $n-p^{\al-1}<\varphi(p^\al)=p^\al-p^{\al-1}$.

Below we let $n\gs p^\al$ and
assume that the desired result holds for smaller values of $n$ not less than
$p^{\al-1}$.
Note that $n'=n-\varphi(p^{\beta})<n$ and also
$n'\gs n-\varphi(p^\al)\gs p^{\al-1}$.

Let $r$ be any integer, and set
$$S=\sum_{k\eq r\,(\mo\ p^{\beta})}
\bi nk(-1)^kf\l(\l\lfloor\f{k-r}{p^{\al}}\r\rfloor\r).\tag3.2$$
By the Chu-Vandermonde identity,
$$\align S=&\sum_{k\eq r\,(\mo\ p^{\beta})}
\sum_{j=0}^{\varphi(p^{\beta})}\bi{\varphi(p^{\beta})}j\bi{n'}{k-j}
(-1)^kf\l(\l\lfloor\f{k-r}{p^{\al}}\r\rfloor\r)
\\=&\sum_{j=0}^{\varphi(p^{\beta})}\bi{\varphi(p^{\beta})}j\sum_{k\eq r\,(\mo\ p^{\beta})}
\bi{n'}{k-j}(-1)^{k}f\l(\l\lfloor\f{k-j-(r-j)}{p^{\al}}\r\rfloor\r)
\\=&\sum_{j=0}^{\varphi(p^{\beta})}\bi{\varphi(p^{\beta})}j(-1)^jS_j,
\endalign$$
where
$$S_j=\sum_{k\eq r-j\,(\mo\ p^{\beta})}
\bi{n'}{k}(-1)^{k}f\l(\l\lfloor\f{k-(r-j)}{p^{\al}}\r\rfloor\r).\tag3.3$$

For any $j=0,1,\ldots,\varphi(p^{\beta})$, by the induction hypothesis we have
$$\ord_p(S_j)\gs \gamma=\l\lfloor\f{n'-p^{\al-1}-l\da_{\al,1}}{\varphi(p^{\al})}\r\rfloor
-(l-1)\al-\beta,$$
and Lemma 2.2 yields that
$$\bi{\varphi(p^{\beta})}j\eq\cases(-1)^j\ (\mo\ p)&\t{if}\ p^{\beta-1}\mid j,
\\0\ (\mo\ p)&\t{if}\ p^{\beta-1}\nmid j.\endcases$$
Thus, if $\gamma\gs0$ then
$$S\eq\sum_{j=0}^{p-1}\bi{\varphi(p^{\beta})}{p^{\beta-1}j}(-1)^{p^{\beta-1}j}S_{p^{\beta-1}j}
\eq\sum_{j=0}^{p-1}S_{p^{\beta-1}j}
\ \l(\mo\ p^{\gamma+1}\r).$$
Observe that
$$\sum_{j=0}^{p-1}S_{p^{\beta-1}j}
=\sum_{k\eq r\, (\mo\ p^{\beta-1})}
\bi{n'}k(-1)^kf\l(\l\lfloor\f{k-(r-p^{\beta-1}j_k)}{p^{\al}}\r\rfloor\r),$$
where $j_k$ is the unique integer in $\{0,\ldots,p-1\}$
with $p^\beta\mid k-(r-p^{\beta-1}j_k)$.
For $k\eq r\ (\mo\ p^{\beta-1})$, clearly
$$\f{k-r+p^{\beta-1}j_k}{p^\beta}
=\f{k-r'-p^{\beta-1}(p-1-j_k)}{p^\beta}=\l\lfloor\f{k-r'}{p^\beta}\r\rfloor$$
where $r'=r-\varphi(p^\beta)$.
Therefore $\sum_{j=0}^{p-1}S_{p^{\beta-1}j}=S',$
where $$S'=\sum_{k\eq r'\, (\mo\ p^{\beta-1})}
\bi{n'}k(-1)^kf\l(\l\lfloor\f{k-r'}{p^{\al}}\r\rfloor\r).\tag3.4$$
From the above it follows that
$$\ord_p(S-S')\gs\gamma+1
\gs\l\lfloor\f{n-p^{\al-1}-l\da_{\al,1}}{\varphi(p^\al)}\r\rfloor-(l-1)\al-\beta.$$

Let $l_0=l$ if $\al=1$, and $l_0=\min\{l,\da_{\beta-1,0}\}$ if $\al>1$.
As $w_l(\al,\beta-1)<w_l(\al,\beta)=w$, by the induction hypothesis we have
$$\align\ord_p(S')\gs&
\l\lfloor\f{n'-p^{\al-1}-l_0}{\varphi(p^{\al})}\r\rfloor-(l-1)\al-(\beta-1)
\\\gs&\l\lfloor\f{n-p^{\al-1}-l\da_{\al,1}}{\varphi(p^\al)}\r\rfloor-(l-1)\al-\beta.
\endalign$$
(Note that if $\al>1=\da_{\beta-1,0}$ then $\beta=1<\al$
and hence $n'-1+\varphi(p^\al)\gs n'+\varphi(p^\beta)=n$.)

Combining the above we finally obtain that
$$\ord_p(S)=\ord_p((S-S')+S')
\gs\l\lfloor\f{n-p^{\al-1}-l\da_{\al,1}}{\varphi(p^\al)}\r\rfloor-(l-1)\al-\beta.$$
Since $\da_{\beta,0}=0$, this concludes the induction step in Case 2.

The proof of Theorem 1.1 is now complete.

\medskip

\Ack. The work was done during the author's visit to the University of California
at Irvine, and he would like to thank Prof. Daqing Wan for the kind invitation.
The author is also indebted to the referee for helpful comments.

\bigskip

\leftline{Department of Mathematics (and Institute of Mathematical Science)}
\leftline{Nanjing University}
\leftline{Nanjing 210093}
\leftline{People's Republic of China}
\leftline {E-mail: zwsun\@nju.edu.cn}
\leftline {{\tt http://pweb.nju.edu.cn/zwsun}}

\enddocument